\documentclass[12pt]{amsart}
\usepackage{amssymb, latexsym, amsmath, amscd, oldlfont}
\usepackage[dvips]{graphicx}
\usepackage{psfrag}

\advance\hoffset by -0.5 truecm

\def\UU{\widetilde{U}}
\def\H{{\cal H}}
\def\HH{{\widetilde{H}}}

\def\ppi^n{\widetilde{\pi}}

\def\ep{{\epsilon}}

\def\om{{\omega}}

\def\vr{\varphi}

\def\ppi{\widetilde{\pi}}

\def\hra{\hookrightarrow}

\def\HZ{\operatorname{HZ}}

\def\R{{\mathbb R}}

\def\T{{\mathbb T}}
\def\CP{{\mathbb C}{\mathbb P}}

\def\S{{\mathbb S}}

\newtheorem{theorem}{Theorem}[section]

\newtheorem{lemma}{Lemma}[section]
\newtheorem{corollary}{Corollary}[section]
\theoremstyle{definition}
\newtheorem{definition}{Definition}[section]

\def\thebibliography#1{\section*{References}\list
 {[\arabic{enumi}]}{\settowidth\labelwidth{[#1]}\leftmargin\labelwidth
 \advance\leftmargin\labelsep
 \usecounter{enumi}}
 \def\newblock{\hskip .11em plus .33em minus .07em}
 \sloppy\clubpenalty4000\widowpenalty4000
 \sfcode`\.=1000\relax}

\newenvironment{remark}{\par\medskip\noindent{\bf 
Remark.}}{\hfill\par\medskip}
\newenvironment{remarks}{\par\medskip\noindent{\bf 
Remarks.}}{\hfill\par\medskip}

\author{Leonardo Macarini}
\title[Hofer-Zehnder capacity of standard cotangent bundles]{Hofer-Zehnder capacity of\\
standard cotangent bundles}
\address{Instituto de Matem\'atica Pura e Aplicada - IMPA\\
         Estrada Dona Castorina, 110 - Jardim Bot\^anico\\
         22460-320 Rio de Janeiro RJ\\
         Brasil}
\email{leonardo@impa.br}         
\date{current version: August 2003}
\thanks{This work was partially CNPq-Profix, Brazil.}

\begin{document}

\begin{abstract}
Let $M$ be a compact manifold with an effective semi-free circle action whose
fixed point set has trivial normal bundle. We prove that its cotangent bundle
endowed with the canonical symplectic form has bounded Hofer-Zehnder sensitive
capacity. We give several examples like the product of any compact manifold
with $\S^n$ or a connected sum $\CP^n \# ... \# \CP^n$.
\end{abstract}

\maketitle

\section{Introduction}

A remarkable fact in the modern development of Symplectic Topology was the introduction by H. Hofer and E. Zehnder of certain symplectic invariants relating the classical problem of the existence of periodic orbits on prescribed energy levels of Hamiltonian systems with symplectic rigidity phenomena \cite{HZ}. It is defined in the following way:

\begin{definition}
\label{defcHZ}
Given a symplectic manifold $(M,\om)$ define the Hofer-Zehnder capacity of $M$
by
$$ c_{\HZ}(M,\om) = \sup_{H \in \H_a(M,\om)}\ \max H, $$
where $\H_a(M,\om)$ is the set of {\it admissible Hamiltonians} $H$ defined on
$M$, that is,
\begin{itemize}
\item $H \in \H(M) \subset C^\infty(M,\R)$, where $\H(M)$ is the set
of {\it pre-admissible Hamiltonians} on $M$, that is, $H$ satisfies the
following properties: $H \geq 0$, there exist an open set 
$V \subset M$ such that $H|_V \equiv 0$ and a compact set 
$K \subset M\setminus\partial M$ satisfying $H|_{M\setminus K}\equiv \max H$;
\item every non-constant periodic orbit of the Hamiltonian vector field $X_H$ has period greater than 1.
\end{itemize}
\end{definition}

It is easy to prove that if a symplectic manifold has {\it bounded Hofer-Zehnder capacity}, that is, if $c_{\HZ}(U,\om) < \infty$ for every open subset $U \subset M$ with compact closure, then given any Hamiltonian $H: M \to \R$ with compact energy levels, there exists a dense subset $\Sigma \subset H(M)$ such that for every $e \in \Sigma$ the energy hypersurface $H^{-1}(e)$ has a periodic solution. In particular, the Weinstein conjecture is true for such manifolds.

Actually, if there exists a neighborhood of an energy hypersurface $H^{-1}(e_0)$ with finite Hofer-Zehnder capacity then there exists $\ep>0$ and a dense subset $A \subset (e_0-\ep,e_0+\ep)$ such that $H^{-1}(e)$ has periodic orbits for every $e \in A$. This result was improved in \cite{MS}, where it is shown that there are periodic orbits on $H^{-1}(e)$ for almost all $e \in (e_0-\ep,e_0+\ep)$ with respect to the Lebesgue measure.

We will consider here a refinement of the original Hofer-Zehnder capacity considering periodic orbits with homotopy class in a prescribed subset $\Gamma$ of $\pi_1(M)$ \cite{Mac1,Mac2,Sch}:

\begin{definition}
Given a symplectic manifold $(M,\om)$ and a subset $\Gamma \subset \pi_1(M)$
define the {\it Hofer-Zehnder $\Gamma$-semicapacity} (or {\it Hofer-Zehnder $\Gamma$-sensitive capacity}) of $M$ by
$$ c_{\HZ}^\Gamma(M,\om) = \sup_{H \in \H_a^\Gamma(M,\om)} \max H, $$
where $\H_a^\Gamma(M,\om)$ is the set of {\it $\Gamma$-admissible Hamiltonians} defined
on $M$, that is,
\begin{itemize}
\item $H \in \H(M)$, that is, $H$ is pre-admissible;
\item every nonconstant periodic orbit of $X_H$ whose homotopy class belongs to
$\Gamma$ has period greater than 1.
\end{itemize}
\end{definition}

Obviously, $c_{\HZ}(M,\om) = c_{\HZ}^{\pi_1(M)}(M,\om) \leq c^\Gamma_{\HZ}(M,\om)$ for every $\Gamma \subset \pi_1(M)$. Moreover, it is easy to see that if the Hofer-Zehnder $\Gamma$-sensitive capacity is bounded then the periodic orbits given by the discussion above have homotopy class in $\Gamma$.

Unfortunately, it is very hard, in general, to prove that a symplectic manifold
has bounded Hofer-Zehnder capacity. For instance, for a particularly important
class of symplectic manifolds, namely, standard cotangent bundles, very little
is known. Actually, the most general class of known examples
are cotangent bundles of manifolds endowed with free circle actions
\cite{Mac2}. It generalizes previous results for product manifolds $M \times
\S^1$ \cite{HV,Lu1,Lu2} and the torus $\T^n$ \cite{HZ,Jia}.

The goal of this note is to give a larger class of manifolds for which the
standard cotangent bundle has bounded Hofer-Zehnder sensitive capacity. Let us
remember that a group action $G \times M \to M$ is {\it semi-free} if it is
free outside the fixed point set and it is {\it effective} if for any $g \neq
e$, there exists $x \in M$ such that $gx \neq x$. 

\begin{theorem}
\label{thmA}
Let $M$ be a compact manifold with an effective semi-free circle action $\vr$. Suppose that its fixed point set has trivial normal bundle. Then its standard cotangent bundle has bounded Hofer-Zehnder $G_\vr$-sensitive capacity, where $G_\vr \subset \pi_1(T^*M) = \pi_1(M)$ is the subgroup generated by the homotopy class of the orbits of $\vr$.
\end{theorem}

\begin{remarks}
\begin{itemize}
\item It is a well know fact that each connected component of the fixed point set of a circle action is a smooth submanifold.
\item Actually, what we need in the proof is the existence of a vector field on $M$ that is nowhere tangent to the fixed point set. It obviously exists if the normal bundle is trivial.
\end{itemize}
\end{remarks}

There are a lot of examples of such manifolds. In particular, we have the following immediate corollaries:

\begin{corollary}[\cite {Mac2}]
Let $M$ be a compact manifold endowed with a free circle action $\vr$. Then its standard cotangent bundle has bounded Hofer-Zehnder $G_\vr$-sensitive capacity.
\end{corollary}

\begin{corollary}
Let $M$ be a product of a compact manifold with the sphere $\S^n$ or a connected sum $\CP^n \# ... \# \CP^n$. Then $T^*M$ has bounded Hofer-Zehnder $0$-sensitive capacity, that is, the periodic orbits are contractible.
\end{corollary}

In fact, it is easy to see that both $\S^n$ and $\CP^n$ have semi-free circle
actions with isolated fixed points. To get this action for a connected sum
$\CP^n \# ... \# \CP^n$, consider a fixed point $p$ in $\CP^n$ and let $U$ be a
neighborhood of $p$ such that we can identify each complex projective space by
a map $\phi: \partial U \to \partial U$ equivariant with respect to the action.
Notice that this argument works for any connected sum $M_1 \# M_2$ such
that both $M_1$ and $M_2$ have a circle action with isolated fixed points.

As an immediate consequence of the previous theorem, we have the following application on the existence of periodic orbits on prescribed energy levels of exact magnetic flows (for an introduction to magnetic flows see, for instance, \cite{CMP,FS,Mac1}):

\begin{corollary}
Let $M$ be a compact Riemannian manifold as in Theorem \ref{thmA} and $\theta$ a 1-form on $M$. Consider the Hamiltonian $H(q,p) = \|p-\theta_q\|^2$ on $T^*M$. Then the Hamiltonian flow of $H$ has periodic orbits with homotopy class in $G_\vr$ on almost all energy levels.
\end{corollary}

\medskip
\noindent
{\bf Acknowledgements.} I am very grateful to Felix Schlenk for very useful comments.

\section{Proof of Theorem \ref{thmA}}

Consider a Riemannian metric on $M$ invariant by the circle action (it exists
because the group is compact) and let $TM$ be the tangent bundle endowed with
the pullback $\om_0$ of the canonical symplectic form by the bundle isomorphism
induced by the metric. Let $\xi$ be the lifted circle action on the tangent
bundle. This action is Hamiltonian and, since the metric is invariant by the
action, it is given by $\xi_{-t}(x,v) = (\vr_t(x),(d\vr_t)_x v)$. Thus, its
fixed point set is the tangent subbundle of the fixed point set $F$ of $\vr$
because the derivative of $\vr$ in the normal direction of $F$ is a rotation
with the same period of the action.

Note that each connected component of the fixed point set has dimension
strictly less than the dimension of $M$ because the action is effective and the
group is compact.

\begin{lemma}
Let $U \subset TM$ be an open subset with compact closure. Then there exists an
open subset with compact closure $\UU \subset TM\setminus TF$ such that given a
$G$-admissible Hamiltonian $H$ defined on $U$, there exists a $G$-admissible
Hamiltonian $\HH$ defined on $\UU$ satisfying $\max \HH = \max H$.
\end{lemma}

\begin{proof}
Since $F$ has trivial normal bundle, there exists an exact 1-form $\theta$ on
$M$ whose dual (gradient) vector field $X$ satisfies $X(x) \in T^\perp_x F$ and
$X(x) \neq 0$ for every $x \in F$. Moreover, since $U$ has compact closure,
there exists a constant $a>0$ sufficiently great such that $TF+aX$ does not
intersect $U$.

\begin{remark}
Note that the existence of $\theta$ already follows from the existence of a vector field $Y$ on $M$ nowhere tangent to $F$. Indeed, consider the submanifold $S$ given by $\cup_{-\ep<t<\ep, x \in F} \exp_x(tY(x))$ for $\ep>0$ sufficiently small and define a smooth function $f$ on $S$ such that $df(Y(x))>0$. Finally, extend smoothly the function $f$ to all of $M$ and define $\theta = df$.
\end{remark} 

We can suppose, without loss of generality, that $H$ is defined on $M$ such that $H|_{M\setminus U} \equiv \max H$. Consider the Hamiltonian $\HH$ given by $\HH(x,v) = H(x,v+aX)$. We have that $\HH$ is also $G$-admissible because the translation by $aX$ along the fiber is a symplectomorphism (since $\theta$ is closed) and homotopic to the identity. Notice that $\HH$ is defined on $\UU:=U - aX$, which, by the choice of $a$, does not intersect $TF$.
\end{proof}

Now, we need the following result that corresponds to a particular case of a more general theorem proved in \cite{Mac2}:

\begin{theorem}[\cite{Mac2}]
\label{action} 
Let $(M,\om)$ be a geometrically bounded exact symplectic manifold and $N \subset M$ an open subset that admits a free Hamiltonian circle action $\vr$. Then, given any open subset $U \overset{i}{\hra} N$ with compact closure,
$$ c_{\HZ}^{i_*^{-1}G_\vr}(U,\om) < \infty, $$
where $i_*: \pi_1(U) \to \pi_1(N)$ is the induced homomorphism on the
fundamental group and $G_\vr \subset \pi_1(N)$ is the subgroup generated by the
orbits of the circle action.
\end{theorem}

We refer the reader to \cite{ALP} for the definition and a discussion of the concept of geometrically bounded symplectic manifold. Here we only mention that it is a well know fact that standard cotangent bundles are geometrically bounded.

Consequently, combining the previous lemma and the theorem above, we conclude that
$$ c_{\HZ}^{G_\vr}(U,\om_0) \leq c_{\HZ}^{G_\vr}(\UU,\om_0) < \infty, $$
as desired.

\end{document}